\documentclass[12pt]{amsart}
\usepackage{amssymb,amsmath,amscd,graphicx,amsthm,color,enumerate}
\usepackage[square, numbers, comma, sort&compress]{natbib}
\usepackage{pst-node}
\usepackage{tikz-cd}

\theoremstyle{definition}
\newtheorem{thm}{Theorem}[section]
\newtheorem{lem}[thm]{Lemma}
\newtheorem{prop}[thm]{Proposition}
\newtheorem{defn}[thm]{Definition}
\newtheorem{cor}[thm]{Corollary}
\newtheorem{conj}[thm]{Conjecture}

\theoremstyle{definition}
\newtheorem*{pf}{Proof}

\theoremstyle{remark}

\newcommand{\p}[1]{[#1]}
\newcommand{\plucker}[4]{\frac{\p{#1}\p{#2}}{\p{#3}\p{#4}}}

\begin{document}
\title{On bounded ratios of minors of totally positive matrices}

\author{DANIEL SOSKIN}
\author{MICHAEL GEKHTMAN}

\begin{abstract}
We construct several examples of bounded Laurent monomials in minors of an $n\times n$ totally positive matrix, which can not be factored into a product of so called primitive bounded ratios thus disproving the conjecture about factorization of bounded ratios stated in \cite{fallat2003multiplicative} by Fallat, Gekhtman, and Johnson. However, all of the found examples still satisfy the subtraction-free property also conjectured in \cite{fallat2003multiplicative}. In addition, we show that the set of all of the bounded ratios forms a polyhedral cone of dimension $\binom{2n}{n}-2n$.
\end{abstract}

%\subjclass[2010]{53D17,13F60}
\keywords{totally positive matrices, determinantal inequalities}

\maketitle
%\tableofcontents

\section{Introduction}
A real $n\times n$ matrix $A$ is called \emph{totally positive} if every minor of $A$ is positive. Initially, total positivity arose in three different areas. It was studied by Gantmacher and Krein in small oscillations of vibrating systems (see \cite{GKOsz},\cite{GKOsc}), by Sch\"onberg in applications to analysis of real roots of polynomials and spline functions (see \cite{cur}), and by Karlin in integral equations and statistics (see \cite{ando87}, \cite{gas2013total}). Totally positive matrices appear in many other areas of mathematics (see, e.g., \cite{FJTNNMatrices}).
%$\quad$ 
For $I,I' \subseteq \{1,2,\ldots,n\}$ with
$|I|=|I'|$, we denote the minor of $A$ with row set $I$ and column
set $I'$ as $\Delta_{I,I'}(A):=\det A(I|I')$.
We follow the convention that
$\Delta_{\emptyset,\emptyset}(A)=1$. We are interested in ratios of products of minors of $A$
\begin{equation}\label{eq:R}R=\frac{\Delta_{I_{1},I'_{1}}(A)\Delta_{I_2,I_2'}(A)...\Delta_{I_p,I_p'}(A)}{\Delta_{J_1,J_1'}(A) \Delta_{J_2,J_2'}(A)...\Delta_ {J_q,J_q'}(A)},
\end{equation}  
more precisely, in conditions on  $R$ to be bounded on the locus of totally positive elements in $GL_{n}$.

%\textcolor{red}{Recall peper then this q was stated}
This question was first posed in \cite{fallat2000det}. In \cite{fallat2003multiplicative} a large class of bounded ratios
of products of principal minors was classified. In particular, necessary
and sufficient conditions were stated for a ratio of products of two minors to
be bounded over totally positive matrices.

This result was generalized to non-principal minors in
\cite{skandera2004inequalities}.
Necessary condition in
\cite{fallat2003multiplicative} was generalized in \cite{boocher2008generators} to the case of non-principal minors and an explicit factorization of ratios of products of two minors
into products of so-called {\em primitive ratios} was constructed.
It was conjectured in \cite{fallat2003multiplicative} that all bounded ratios can be factored into products
of primitive ratios
\cite{fallat2003multiplicative}.

One of the goals of this paper is to disprove this conjecture. Counterexamples are presented in Theorem~\ref{t:main}. To this end we study the cone of all bounded ratios. In particular, in Theorem~\ref{t:polyh} we show that this cone is polyhedral. The main tool we employ is a \emph{network parameterization} of totally  positive matrices. 

The relationship between totally positive matrices and \emph{directed
acyclic weighted planar networks} is well studied.  It was first discussed by Karlin and
McGregor in 1959 \cite{karlin1959coincidence} and later developed by Gessel and Viennot \cite{Gessel-Viennot} who built on the work of Lindstr\"om \cite{Lindstrom}. For a more modern presentation, we refer to \cite{fomin2000total} and references therein. Postnikov introduced transition from directed acyclic weighted planar networks to the \emph{planar bicolored graphs} with assigned \emph{face weights}
\cite{postnikov2006total}.
Lately, the theory of plabic graphs has attracted a lot of interest in several areas of mathematics and theoretical physics. A good exposition of the subject can be found in \cite{fomin2021introduction}. Below, we use planar networks to better understand the cone $\mathcal{C}_{n}$ of all bounded ratios and use the gained insights to develop computational tools that allow us to provide examples of extremal rays of $\mathcal{C}_{n}$ that correspond to baounded ratios disproving the conjecture of Fallat-Johnson-Gekhtman.

The paper is organized as follows.
Section 2 contains the necessary background information and two conjectures on bounded ratios which are discussed later on.
Section 3 contains preliminary technical results. We show that any unbounded ratio can be detected by some family of face weights of type $t^{\gamma_{i}}$. Then we show that the vectors of exponents of the bounded ratios form a polyhedral cone in a vector space $V_n$ formally spanned by Pl\"ucker coordinates for the Grassmannian $Gr(n,2n)$. This cone is given by a system of linear homogeneous inequalities. For matrices of order $4$ the system was computed in Matlab using the {\tt qskeleton} software \cite{bastrakov2015fast}. 
In Section 4 we study properties of the primitive vectors. In particular, we describe simple linear relations among them. These relations allow us to find dual description of the cone generated in $V_n$ by {\em primitive  vectors}, i.e. vectors that correspond to primitive ratios. It is given by a system of linear homogeneous inequalities (even on a supercomputer, {\tt qskeleton}  was not able to solve the system in a reasonable time). By comparing it with the system of inequalities corresponding to bounded ratios, we develop computational methods for searching for new bounded ratios.

In Section 5 we describe two numerical methods which result in a list of new extreme rays/bounded ratios. 

\section{Preliminaries}
In this section we review some background information about planar networks, Grassmannians, and bounded ratios of minors of totally positive matrices. 
\subsection{Planar networks with face weights and totally positive matrices}
We will consider a particular example of a directed acyclic planar network with weighted faces \cite{postnikov2006total}. %\cite{fomin1999total}. 
It is a planar directed graph with $n$ labeled sources on the left and
$n$ labeled sinks
at the bottom. Each edge is oriented either to the right or down. Each face contains a positive weight (see Fig.1).

Network of size $n$ has $n^2$ face weights indexed as in Fig.1. We denote it by $\mathcal{N}_{n}$. Networks $\mathcal{N}_{n}$ can be used to parameterize the set $M^{>0}_{n}$ of square totally positive matrices of order $n$ by associating to each weighted network
%Both sources and sinks are labeled bottom to top.
%For each edge with a red dot (see Fig.1) as2 positive real weight is associated. All edges without dots have weights 1.
%\begin{figure}[htb]
%\centering
%\includegraphics[width=4%.2in]{general_ntw3.png}
%\caption{\label{fig:general}General Planar Network}
%\end{figure}
a matrix $X$ in the following way:  \begin{equation}x_{ij}=\sum_{\pi:i\to
j} w(\pi)\label{bijection} \end{equation} where the summation is
over all directed paths $\pi$ connecting the source $i$ and the sink $j$. Here $w(\pi)$ is product of all weights assigned to faces located to the right of the path $\pi$. The following theorem is a particular case of a general result in \cite{postnikov2006total}.
\begin{thm}
The set of planar networks $\mathcal{N}_{n}$ with all possible positive weights is in bijection with the set of totally positive square matrices of order $n$.
\end{thm}
The theorem above relies on the formula for the minors of $X$ in terms of positive weights.
\begin{lem}[Lindstr\"om’s Lemma \cite{Lindstrom,Gessel-Viennot}]
Any minor $\Delta_{I,J}$
of weighted matrix $X$ of a planar network is given by the formula:
\begin{equation}\Delta_{I,J}=\sum_{P_{I\to
J}} w(P_{I\rightarrow J})\label{Lind} \end{equation} where $P_{I \rightarrow J}$ is a family of non-intersecting paths and $w(P_{I\rightarrow J})$ is weight of the path family obtained as a product of weights of all paths in the family.
\end{lem}

\begin{figure}
\centering
\includegraphics[width=3.in]{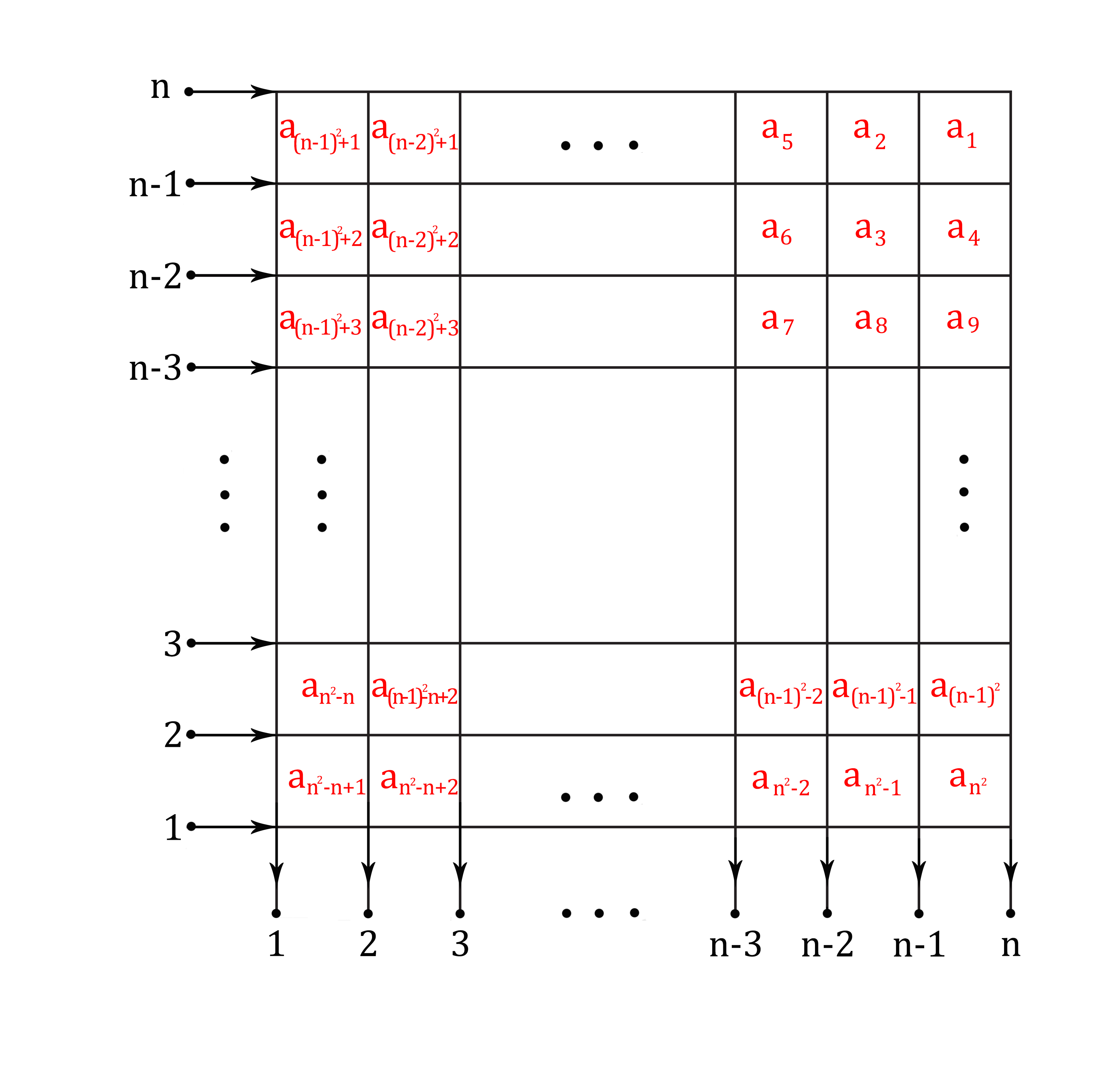}
\label{fig:f1}
\caption{
%Directed planar network with face weights
\\}
\end{figure}
\subsection{Grassmannians}
The {\em real Grassmannian} $Gr(n,2n)$ is the manifold
of $n$-dimensional subspaces in $\mathbb {R}^{2n}$. An element $P \in Gr(n,2n)$ is represented by a $2n \times n$-matrix $B$ whose columns span $P$. The \emph{Pl\"ucker coordinates of $P$ with respect to $B$} are defined as a vector of all $n\times n$ minors of the matrix $B$ viewed as an element of the real projective $({2n \choose n}-1)$-space. 

The \emph{totally positive Grassmannian} $Gr^{>0}(n,2n)$ is a subset of $Gr(n,2n)$ consisting of all points that have a matrix representative with all positive Pl\"ucker
coordinates \cite{postnikov2006total}.

Let $[n]$ be an interval of natural numbers from $1$ to $n$ (more generally, $[a,b]$ is an interval of all natural numbers $n:a\leq n \leq b$ ).

We consider an embedding of $M_{n}^{>0}$ into $Gr^{>0}(n,2n)$. A matrix
$X\in M_{n}^{>0}$ corresponds to the point in Grassmannian with a representative  
$$\overline{X}= \begin{bmatrix}
X\\
W_{0}
\end{bmatrix},~\text{where}~ W_{0}=\Big((-1)^{i+1}\cdot\delta_{j,n-i+1}\Big)^{n}_{i,j=1}~~~~~\text{i.e.}~~w_{ij}=\begin{cases}
      -1^{i+1} & \text{if $j=n-i+1$}, \\
      0&\text{otherwise.}
    \end{cases}$$
%\\ where $w_{ij}=-1^{i+1}$ if $j=n-i+1$, and $0$ otherwise. \textcolor{red}{write as a system}
Then $\Delta_{I,I'}(X)=\Delta_{I",[n]}(\overline{X})=:\Delta_{I"}(\overline{X})$, 
where $I"={I\cup \{2n+1-i|i \in [n] \setminus I'\}}$. 
Observe that $\Delta_{I"}(\overline X)$ $=\Delta_{(I" \cap [n]),([n+1,...,2n]\setminus I")}(X)$ for any set $I"$ (\cite{skandera2004inequalities}). Due to this correspondence, instead of ratios of minors (\ref{eq:R}) we can study ratios of Pl\"ucker coordinates.
\subsection{Background information about primitive ratios and bounded ratios}
The simplest bounded ratios were defined in \cite{boocher2008generators} and were called \emph{`basic` ratios}. We call them \emph{`primitive` ratios} to avoid possible confusion. 

\begin{defn}[\cite{boocher2008generators}]\label{defn:pr}A primitive ratio is one of the form
\[ \plucker{i,j+1,\Delta}{i+1,j,\Delta}{i,j,\Delta}{i+1,j+1,\Delta} \]
where $i,j\in \{1,\ldots,2n\}$ and $\Delta\subset \{1,\ldots,2n\}$
such that $|\Delta|=n-2$ and $i$, $i+1$, $j$, $j+1$ and $\Delta$ are
all distinct, $i<j$.  Here index $j+1$ is understood mod $2n$. Denote this ratio by $R_{i,j,\Delta}$.
\end{defn}

Let $\alpha$ be a sequence of index sets in the numerator and $\beta$ be a sequence of index sets in the denominator of a ratio of Plucker coordinates.
\begin{defn} [\cite{boocher2008generators}]
For $i\in\{1,2,\ldots,2n\}$, let
$f_\alpha(i)$ be the number of index sets in $\alpha$ that contain
$i$.  
%\begin{defn}[ST0 Condition]
If
$f_\alpha(i)=f_\beta(i)$ for all $i$, we say the ratio
$\alpha/\beta$ satisfies the ST0 condition.
\end{defn}

\begin{lem}[\cite{boocher2008generators}]
If a ratio is bounded for all totally positive matrices, then the ratio satisfies the ST0
condition. 
\end{lem}

\begin{thm}
\label{thm:symm}[\cite{boocher2008generators},\cite{fallat2003multiplicative}] For an index set $\alpha_j$, we define a cyclic shift
$\sigma(\alpha_j) = \{ i+1 \mod{2n}\, | \, i\in \alpha_j\}$ and a reflection $\rho(\alpha_j) = \{ (2n+1) -i \, | \, i \in \alpha_j \}$. For a sequence $\alpha=(\alpha_1,\ldots,\alpha_p)$ of index sets,
we define $\sigma(\alpha)=(\sigma(\alpha_1),\ldots,\sigma(\alpha_p))$ and define $\rho(\alpha)=(\rho(\alpha_1),\ldots,\rho(\alpha_p)).$ Then 
\\ (1) the ratio $\alpha/\beta$ is bounded if
and only if $\sigma(\alpha)/\sigma(\beta)$ is bounded;
\\ (2) the ratio $\alpha/\beta$ is bounded if
and only if $\rho(\alpha)/\rho(\beta)$ is bounded.
\end{thm}
Any ratio of minors 
%$\alpha/ \beta$ 
is a rational function $\frac{p}{q}$ in terms of the face weights of the network parametrization. This leads to the following
\begin{defn}
A ratio  is called \emph{subtraction-free} if $q-p$ is a polynomial function in face weights with all positive coefficients.
\end{defn}

\begin{prop}[\cite{skandera2004inequalities},\cite{boocher2008generators}]
Every bounded ratio of minors of a $3\times 3$ totally positive matrix can be
written as a product of positive powers of the primitive ratios.
Furthermore, every such bounded ratio is substraction-free and bounded by
1.
\end{prop}

For matrices of order $n$ there are $N=\binom{2n}{n}$ Pl\"ucker coordinates, for which we chose a lexicographical order: $\Delta_{1}$ is $\Delta_{[1,2,3...n]}$,...,$\Delta_{N}$ is $\Delta_{[n+1,n+2...2n]}$. Bounded ratios (\ref{eq:R}) can be viewed as bounded Laurent monomials
\begin{equation}\label{eq:Rlor}
R=\Delta_{1}^{{\alpha_{1}}}\cdot...\cdot\Delta_{N}^{{\alpha_{N}}},\quad \alpha_{i} \in \mathbb{Z} \end{equation} 
Let us consider a linear space $V_n=\mathbb {R}^{N}$, where each axis $x_{i}$ corresponds to $\Delta_{i}$. Let $e_{I}$ denote the basis of $V_n$ indexed by $n$-element subsets $I \in \binom{[2n]}{n}$. Then, each monomial $(\ref{eq:Rlor})$ corresponds to a vector $\Vec{v}=(\alpha_{1},...,\alpha_{N})$ in $V_n$. Denote the vector associated with $R_{i,j,\Delta}$ by $\Vec{v}_{i,j,\Delta}$.
Denote the set of vectors $\Vec{v}_{i,j,\Delta}$ with $|\Delta|=n-2$ by $\mathcal{P}_{n}$.
Below, depending on the context, we will use both monomial and vector descriptions.

\subsection{Conjectures on bounded ratios}
\begin{conj}[Fallat, Gekhtman, Johnson \cite{fallat2003multiplicative}]\label{conj:mainconj}
Let $R$ be a ratio of minors, then
\\(1) $R$ is bounded if and only if $R$ can be written as a
product of primitive ratios.
\\(2) $R$ is bounded if and only if $R$ is subtraction
free.
\end{conj}
One of the main goals of this paper is to disprove part (1) of this conjecture. 

\section{Description of the set of the bounded ratios}
We denote the set of vectors in $V_n$ that corresponds to the set of all bounded ratios for $M_{n}^{>0}$ by $\mathcal{C}_{n}$. Clearly, $\mathcal{C}_{n}$ is a cone.
In this section, we show that in order to determine $\mathcal{C}_{n}$, it is enough to consider face weights of type $t^{\lambda_{i}}$, where  $\lambda_{i} \in \mathbb{Z}$ are constants and $t \in \mathbb{R}^{+}$. This will allow us to prove that $\mathcal{C}_{n}$ is a polyhedral cone in the vector space $V_n$. 

\subsection{Specialization of weight parameters} 

\begin{defn} Let $f(x_{1},... ,x_{m})=\sum_{k}c_{k}x_{1}^{\gamma_{1}(k)}\cdots x_{m}^{\gamma_{m}(k)}$ be a polynomial. The \emph{Newton polytope} $\mathbf{N}(f)$ of $f$ is a convex hull of all distinct vectors ${\Vec{\gamma}(k)}$. 
%Denote it as $\mathbf{N}(f)$.
\end{defn}
\begin{prop}\label{pr:polyN}
Let $f_{1}(x_1, ... ,x_m)$, $f_{2}(x_1, ... ,x_m)$ be two  polynomials with positive coefficients. 
%Let $\frac{f_{1}}{f_{2}}$ be a function defined for $x_i \in \mathbb{R}^{+}, 1\leq i\leq m$. 
Then $\frac{f_{1}}{f_{2}}$ is bounded on $(\mathbb{R}^{+})^m$ if and only if $\mathbf{N}(f_{1})\subseteq\mathbf{N}(f_{2})$.  
\end{prop}
\begin{proof}
Let $f_{2}=\sum_{k}c_{k}x_{1}^{\gamma_{1}(k)}\cdots x_{m}^{\gamma_{m}(k)}$.
First, we show that if $\mathbf{N}(f_{1})\nsubseteq\mathbf{N}(f_{2})$, then $\frac{f_{1}}{f_{2}}$ is unbounded. Indeed, let $O'$ be a vertex of $\mathbf{N}(f_{1})$ such that $O'\notin \mathbf{N}(f_{2})$. Then there is a hyperplane $\mathbf{h}$ which separates $O'$ and $\mathbf{N}(f_{2})$. Pick a point $P \in \mathbf{h}$ such that line $O'P$ is perpendicular to $\mathbf{h}$. For a vector $\Vec{p}=\overline{PO'}=(p_1,... ,p_m)$, consider a specialization $x_{i}=t^{p_i}, 1\leq i\leq m$. Then, any monomial $x_{1}^{\gamma_{1}}... ~x_{m}^{\gamma_{m}}$ is equal to $t^{\gamma_{1}p_{1}}... ~t^{\gamma_{m}p_{m}}=t^{\gamma_{1}p_{1}+...~+\gamma_{m}p_{m}}=t^{\Vec{\gamma}\cdot \Vec{p}}$. We observe that $\overline{OO'}\cdot \Vec{p}=\overline{OP}\cdot \Vec{p}+\overline{PO'}\cdot \Vec{p}=\overline{OP}\cdot \Vec{p}+\Vec{p}\cdot \Vec{p} > \overline{OP}\cdot \Vec{p}$. On the other hand, if $\Vec{\gamma}(T)$ corresponds to a vertex $T$ of $\mathbf{N}(f_{2})$ then $\Vec{\gamma}(T)\cdot\Vec{p}=\overline{OP}\cdot \Vec{p}+\overline{PT}\cdot \Vec{p}< \overline{OP}\cdot \Vec{p}$. Since $\lim_{t\rightarrow \infty}\frac{t^{\overline{OO'}\cdot \Vec{p}}}{t^{\Vec{\gamma}(T)\cdot\Vec{p}}}=\infty$ for any $T$ it follows that $\frac{f_1}{f_2}$ is unbounded.

Now we show that if $\mathbf{N}(f_{1})\subseteq\mathbf{N}(f_{2})$ then the function $\frac{f_1}{f_2}$ is bounded. Let $M$ be the number of vertices in $\mathbf{N}(f_{2})$. For any vertex of $\mathbf{N}(f_{1})$ there are $M$ coefficients $(t_{1},...~t_{M}), ~t_{i}\geq0, ~\sum_{i}t_{i}=1$, such that the monomial corresponding to this vertex is 
$$\prod^{m}_{i=1} x_{i}^{t_{1}\gamma_{i}(1)+t_{2}\gamma_{i}(2)+~...~+t_{M}\gamma_{i}(M)}$$
Here $\Vec\gamma(k)=(\gamma_{1}(k),\ldots, \gamma_{m}(k))$ is the $k$th vertex of $\mathbf{N}(f_{2})$.
It is sufficient to show that the following inequality holds
\begin{equation}\label{eq:conv}
\prod^{m}_{i=1} x_{i}^{t_{1}\gamma_{i}(1)+t_{2}\gamma_{i}(2)+~...~+t_{M}\gamma_{i}(M)} \leq \left(\sum^{M}_{k=1}x_{1}^{\gamma_{1}(k)}\cdots x_{m}^{\gamma_{m}(k)}\right)
\end{equation} 
Indeed,
$$
\mbox{ln}\left(\prod^{m}_{i=1} x_{i}^{t_{1}\gamma_{i}(1)+t_{2}\gamma_{i}(2)+~...~+t_{M}\gamma_{i}(M)}\right)=\mbox{ln}\left(\prod^{M}_{k=1}\left(x_{1}^{\gamma_{1}(k)}\cdots x_{m}^{\gamma_{m}(k)}\right)^{t_k}\right)=
$$
$$=\sum^{M}_{k=1}t_{k}\mbox{ln}\left(x_{1}^{\gamma_{1}(k)}\cdots x_{m}^{\gamma_{m}(k)}\right)\leq \mbox{ln}\left(\sum^{M}_{k=1}t_{k} x_{1}^{\gamma_{1}(k)}\cdots x_{m}^{\gamma_{m}(k)}\right)$$ 
by Jensen's inequality. Since $\mbox{max}_{k}\{t_{k}\}\leq 1$ we have
$$\mbox{ln}\left(\sum^{M}_{k=1}t_{k}x_{1}^{\gamma_{1}(k)}\cdots x_{m}^{\gamma_{m}(k)}\right) \leq \mbox{ln}\left(\sum^{M}_{k=1}x_{1}^{\gamma_{1}(k)}\cdots x_{m}^{\gamma_{m}(k)}\right)$$
The last inequality implies (\ref{eq:conv}).
\end{proof}
\begin{cor}\label{cor:par}
A ratio $R$ is bounded if and only if it is bounded for all face weights $a_{i}$ of the form  $t^{\gamma_{i}}, t>0, \gamma_{i} \in \mathbb{Q}$.
\end{cor}
\begin{proof}
We can apply Proposition~\ref{pr:polyN} to a ratio $R$ expressed in terms of positive face weights. Thus, if $R$ is unbounded there is a parametrization of face weights $a_{i}=t^{\gamma_{i}}$ so that $R$ is arbitrary large for $t$ big enough. Note that in the proof of Proposition~\ref{pr:polyN} we can approximate the values of $\gamma_{i}$ by some rational values, so that all inequalities remain valid.
\end{proof}
%\begin{cor}
%$\mathcal{C}_{n}$ is a cone in the space $V_n$.
%\end{cor}
%\begin{proof}
%Recall  that $R$ can be expressed as follows (4): 
%$$
%R=\Delta_{1}^{{\alpha_{1}}}\cdot...\cdot\Delta_{N}^{{\alpha_{N}}},\quad \alpha_{i} \in \mathbb{Z} 
%$$
%Consider a parametrization of the face weights $a_{i}=t^{\gamma_{i}}$. Then each  $\Delta_{j}$ is a sum of powers of $t$. In particular, the leading monomial of $\Delta_{j}$ is $t^{c_{j}}$, where $c_{j}$ is a constant which depends on $\Vec{\gamma}=(\gamma_1,\ldots,\gamma_{(n-1)^2})$. Then, a linear inequality on $\alpha$'s 
%\begin{equation}
%    I(\Vec{\gamma}): \alpha_{1}c_{1}~+~\alpha_{2}c_{2}~+~...~+\alpha_{N}c_{N}\leq 0
%\end{equation}
%is a necessary condition for $R$ to be bounded. Indeed, if (8) does not hold, then $R$ is arbitrary large for the given parametrization and for a large enough $t$. Due to Corolary~3.3, $R$ is bounded if and only if $R$ satisfies $I(\Vec{\gamma})$ for all $\Vec{\gamma}$. Thus, $\mathcal{C}_{n}$ form a cone in $V_n$.
%\end{proof}
\subsection{The cone of bounded ratios is polyhedral}
%In the following theorem we show that only finitely many inequalities $I(\Vec{\gamma})$ define the cone of bounded ratios, the rest of inequalities are redundant.

%Recall that a polyhedral cone in a linear space is given by a finite system of homogeneous linear inequalities, i.e. it is a cone hull of finitely many extreme rays.
Recall 
%that linear space $V=\mathbb {R}^{N}$, where each axis $x_{i}$ corresponds to $\Delta_{i}$ and 
each ratio $R$ corresponds to a vector $\Vec{v}=(\alpha_{1},...,\alpha_{N})$ in $V_n=\mathbb {R}^{N}$.
\begin{thm}\label{t:polyh}
For any $n$ the cone $\mathcal{C}_{n}$ of all bounded ratios $$R=\frac{\Delta_{I_{1},I'_{1}}(A)\Delta_{I_2,I_2'}(A)...\Delta_{I_p,I_p'}(A)}{\Delta_{J_1,J_1'}(A) \Delta_{J_2,J_2'}(A)...\Delta_ {J_q,J_q'}(A)}$$ is a polyhedral cone in $V_n$ given by a system $F$ of finitely many homogeneous linear inequalities in $\alpha$'s.
\end{thm}

\begin{proof}
Recall  that $R$ can be expressed as follows (\ref{eq:Rlor}): 
$$
R=\Delta_{1}^{{\alpha_{1}}}\cdot...\cdot\Delta_{N}^{{\alpha_{N}}},\quad \alpha_{i} \in \mathbb{Z} 
$$
Evaluate $R$ on an element parameterized by face weights $a_{i}=t^{\gamma_{i}}$. Then each  $\Delta_{j}$ is a sum of powers of $t$ of the form $t^{\Vec{\gamma}\cdot \Vec{u}}$, where $\Vec{u}$ depends on a given path family in the network parametrization and does not depend on $\Vec{\gamma}$. In particular, the leading monomial of $\Delta_{j}$ is $t^{c_{j}}$, where $c_{j}$ is a constant which depends on $\Vec{\gamma}=(\gamma_1,\ldots,\gamma_{n^2})$. Then, a linear inequality on $\alpha$'s 
\begin{equation}\label{eq:inq}
    I(\Vec{\gamma}): \alpha_{1}c_{1}~+~\alpha_{2}c_{2}~+~...~+\alpha_{N}c_{N}\leq 0
\end{equation}
is a necessary condition for $R$ to be bounded. Indeed, if (\ref{eq:inq}) does not hold, then $R$ is arbitrary large for the given parametrization and for a large enough $t$. Due to Corolary~\ref{cor:par}, $R$ is bounded if and only if $R$ satisfies $I(\Vec{\gamma})$ for all $\Vec{\gamma}$.

%We consider parametrization $t^{\gamma_{1}}$,..., $t^{\gamma_{d}}$, where $\gamma_{i}$ are any fixed rational parameters.
%Then, each Pl\"ucker coordinate is a linear combination of powers of $t$ where each power is a linear combination of parameters $\gamma_{i}$. 
For each Pl\"ucker coordinate and each of the corresponding monomials $t^{\Vec{\gamma}\cdot \Vec{u}}$, we have a system of linear inequalities on parameters $\gamma_{i}$ which is necessary and sufficient for  this monomial to be leading for this Plücker coordinate. 
A combination of choices of leading monomials for all Pl\"ucker coordinates is called realizable if corresponding sets of solutions to $N$ systems of linear inequalities have a nonempty intersection. Let us consider a system $S$ which is a union of all $N$ systems for a fixed combination. 
%Then $S$ has a nonempty set of solutions iff the combination of main monomials is realisable. 
We observe that any solution $\Vec{\gamma}$ of $S$ gives exactly one linear inequality $I(\Vec{\gamma})$ on $\alpha_{1}$,...,$\alpha_{N}$, which all bounded $R$ must satisfy. Also, if a solution of $S$ is a linear combination with non-negative coefficients of other  solutions of the same system $S$, then the obtained linear inequality on $\alpha_{i}$ is a linear combination of corresponding inequalities $I(\Vec{\gamma})$ with the same non-negative coefficients. For any system $S$,  the set of solutions is a cone hull of finitely many extreme rays. 
%Note, that extreme rays of the system $S$ are vectors with in the space of parameters $\gamma_{i}$. Since for any realisable combination its corresponding system $S$ has finitely many linear inequalities, the solution set of $S$ has finitely many extreme rays. 
$R$ satisfies all inequalities on $\alpha_{i}$ generated by all solutions to $S$ if and only if $R$ satisfies all inequalities on $\alpha_{i}$ only generated by the union of all extreme rays of all realisable systems $S$, which is a finite set of linear inequalities of type $I(\Vec{\gamma})$.
\end{proof}
\begin{cor}
All extreme rays of $\mathcal{C}_{n}$ have integer coordinates.
\end{cor}
\begin{proof}
We observe that system $S$ of inequalities is a system with integer coefficients. Indeed, $S$ is a union of $N$ systems, where each system defines a leading monomial for a Pl\"ucker coordinate and consists of linear inequalities with integer coefficients. Thus, extreme rays of the cone of solutions of $S$ have rational coordinates. Since there are finitely many of them, they can be scaled in a way that all of them are with integer coordinates. In other words, we consider finite set of parametrizations of type $t^{\gamma_{i}}$ with integer $\gamma_{i}$. Thus, evaluated inequalities of type $I(\Vec{\gamma})$ are inequalities with integer coefficients as well. Thus, the finite set of extreme rays of $\mathcal{C}_{n}$ have rational coordinates, and could be scaled such that all of their coordinates are integer. 
\end{proof}

For each Plücker coordinate for each its monomial we can find explicitly the system on parameters $\gamma_{i}$ which defines whether the monomial is leading. This system can be obtained from straightforward consideration of all possible path-families on corresponding planar network $\mathcal{N}_{n}$ as in Fig.1. Thus, using Matlab and software qskeleton \cite{bastrakov2015fast}, we can find all extreme rays of parameters for all systems $S$ with non-empty sets of solutions. Then, we take the union of all linear inequalities on $\alpha_{i}$'s (which define half-spaces in $V_n$) evaluated on all extreme rays of parameters to form one system $F$. The set of solutions of $F$ is exactly the set of all bounded $R$. 

For $n=4$, we computed the system $F$ explicitly. It is a system of 360 homogeneous linear inequalities in 70 variables. It would be easy to study properties of all bounded ratios if we knew all extreme rays of $F$. Software {\tt qskeleton} \cite{bastrakov2015fast} gives an answer to this type of question, but for smaller systems. For a system as large as $F$, it is practically impossible to use {\tt qskeleton}.

\section{Properties of the primitive ratios}
We denote the cone in $V_n$ generated by all primitive vectors $\Vec{v}_{i,j,\Delta}$ that correspond to primitive ratios (see Definition~\ref{defn:pr}) for $M_{n}^{>0}$ by $\overline{\rm \mathcal{C}}_{n}$.
In order to address Conjecture~\ref{conj:mainconj}, we need to find reduced system of linear inequalities which defines the $\overline{\rm \mathcal{C}}_{n}$ and to compare it to 
%the system of inequalities which defines $\overline{\rm \mathcal{C}}_{n}$ and 
the system of inequalities that defines $\mathcal{C}_{n}$. 
%Note, that the system which defines $\mathcal{C}_{4}$ was already found computationally in Section 3. 

%In order to find the system defining $\overline{\rm \mathcal{C}}_{4}$ we could try to use qskeleton \cite{bastrakov2015fast}, but it is practically impossible because of the size of the input: in case $n=4$ there are 120 primitive  vectors in the 70 dimensional space $V_n$. Fortunately, because of some particular properties of primitive  vectors it is possible to find the system of inequalities defining $\overline{\rm \mathcal{C}}_{4}$ through a different approach. 

\subsection{
%Linear relations on primitive  vectors and 
The  dimension of $\overline{\rm \mathcal{C}}_{n}$}
%In this subsection we study linear dependencies between primitive vectors and show that dimensions of $\overline{\rm \mathcal{C}}_{n}$ and $\mathcal{C}_{n}$ are equal to $\binom{2n}{n}-2n$ for any $n$.

% \begin{figure}[htb]
% \centering
% \includegraphics[width=5.2in]{f4.jpg}
% \caption{\\}
% \label{f4}
% \end{figure}
In order to find  $\mbox{dim}(\mathcal{C}_{n})$ and $\mbox{dim}(\overline{\mathcal{C}}_{n})$ we start with an upper bound.
\begin{lem} Dimension of the cone of the bounded ratios satisfies 
\begin{equation}\label{eq:upb}
\mbox{dim}(\mathcal{C}_{n})\leq\binom{2n}{n}-2n.
\end{equation}
\end{lem}
\begin{proof}
%Upper bound $\mbox{dim}(\mathcal{C}_{n})\leq\binom{2n}{n}-2n$ follows from the ST0 condition (this upper bound also follows from Proposition 4.3).

For each $i\in[2n]$ the condition that a ratio satisfies $f_\alpha(i)=f_\beta(i)$ in ST0 means that the vector corresponding to this ratio belongs to a hyperplane  of codimension 1 in $V_n$. Thus, if a ratio satisfies condition ST0, then the corresponding vector satisfies all $2n$ linear equations. To show that these $2n$ linear equations are independent, let us consider a ratio
$$F_{1}=\frac{[1,2,...,n-1,n][1,2,...,n-1,n+1]...[1,2,...,n-1,2n]}{[1,2,...,n-1,n]^{n+1}[n+1,n+2,...,2n-1,2n]}$$
Then $F_{1}$ satisfies $f_\alpha(i)=f_\beta(i)$ for all $i$ except $i=n$. Thus, the condition $f_\alpha(n)=f_\beta(n)$ is independent from the rest of  conditions $f_\alpha(i)=f_\beta(i)$, where $i$ is not equal to $n$.
Cyclic rotations of $F_{1}$ give independence of the rest of equations. 

Bounded ratios must satisfy ST0, thus $\mathcal{C}_{n}$ is a subset of the intersection of the corresponding hyperplanes. Thus $\mbox{dim}(\mathcal{C}_{n})\leq\binom{2n}{n}-2n$. 
\end{proof}
%Note that $F_{1}$ has degree of the denominator one higher than of the numerator. This is the only occasion, when a considered ratio has different degrees at its numerator and denominator. It was done only for the purpose of showing that all $2n$ equations coming from ST0 are independent. Condition that degrees of numerator and denominator coincide for bounded ratios follows from ST0.

%We can find the dimensions of the $\overline{\rm \mathcal{C}}_{4}$ and $\mathcal{C}_{4}$ just by computing ranks, but it is possible to prove more general statement, that the dimensions of $\overline{\rm \mathcal{C}}_{n}$ and $\mathcal{C}_{n}$ coincide for any $n$.   

\begin{thm}\label{th:dim}
%Dimension of the cones $\overline{\rm \mathcal{C}}_{n}$ and $\mathcal{C}_{n}$ are equal
$$\mbox{dim}(\overline{\rm \mathcal{C}}_{n})=\mbox{dim}(\mathcal{C}_{n})=\binom{2n}{n}-2n.$$
\end{thm}
\begin{proof} 
Since $\overline{\rm \mathcal{C}}_{n}\subseteq\mathcal{C}_{n}$ and due to ~(\ref{eq:upb}) it suffices to show that \begin{equation}\label{eq:ineq}\mbox{dim}(\overline{\rm \mathcal{C}}_{n})\geq\binom{2n}{n}-2n
\end{equation}

Recall that the vector space $V_n$ has a basis $e_{I}$ indexed by $n$-element subsets $I \in \binom{[2n]}{n}$ that correspond to Pl\"ucker coordinates in $Gr(n,2n)$. 

Denote by $\overline{V}_n \subset V_{n}$ the linear span of $\overline{\rm \mathcal{C}}_{n}$. Let $V^{0}_n$ be a $2n$-dimensional subspace of $V_n$ spanned by $e_{I}$ for all $I$ in the set 
$$\mathcal{I}_n = \left \{\{i\} \cup [n+2...2n], (i=1,...,n+1) \right \}\cup\left\{ [n...2n] \setminus \{j\}, (j=n+2,...,2n) \right\}$$
Then (\ref{eq:ineq}) would follow from the following claim: for any $I \in \mathcal{I}^{c}_n = \binom{[2n]}{n} \setminus \mathcal{I}_n$ there exists a vector $w_{I} \in \overline{V}_n$ such that $w_{I}= e_{I}~ (\text{mod}~V^{0}_n)$.

We will prove the claim by induction. For $n=2$, $\mathcal{I}_2 = \{[1,4], [2,4], [3,4], [2,3]\}$, $\mathcal{I}^c_2 = \{[1,2], [1,3]\}$ and $w_{[1,2]}$, $w_{[1,3]}$ are the vectors corresponding to the two primitive ratios $$\frac{\textcolor{red}{[1,2]}[3,4]}{[1,3][2,4]}~~\text{and}~~\frac{[1,4][2,3]}{\textcolor{red}{[1,3]}[2,4]}$$
To execute an induction step, fix $I \in \mathcal{I}^{c}_n$ and pick $p \in [n...2n] \cap I$ and $q \in [1...n+1] \setminus I$ (such $p$ and $q$ clearly exist). Then $I=\{p\} \cup I'$ where $I'$ is an $(n-1)$-element subset in a $2(n-1)$-element set $[1,2n] \setminus \{p, q\}$. The latter can be naturally, identified with $[1,2n-2]$ via an appropriate shift of indices that we denote $\mathbf{S}_{p,q}$. $\mathbf{S}_{p,q}$ induces a natural linear isomorphism between $V^{p,q}_{n-1}:=\text{span}\{e_{I'}:I' \in [1,2n] \setminus \{p,q\}, |I'| = n-1 \}$ and $V_{n-1}$. We retain the notation $\mathbf{S}_{p,q}$ for this isomorphism. If $I' \notin \mathbf{S}^{-1}_{p,q}(\mathcal{I}_{n-1})$ we can apply the induction hypothesis: in $\mathbf{S}^{-1}_{p,q}(\overline{V}_{n-1})$ there exists $\tilde{w}_{I'} = e_{I'}~(\text{mod}~\mathbf{S}^{-1}_{p,q}(V^{0}_{n-1}))$. By restoring the index $p$ to all $(n-1)$ index sets in the equation above, we obtain a vector $\tilde{w}_{I}$ in $\overline{V}_{n}$ such that $\tilde{w}_{I} = e_{I}$ (mod$\tilde{V}^{0}_{n}$) where $\tilde{V}^{0}_{n}$ is spanned by vectors $e_{J}$ with $J={p}\cup\mathbf{S}^{-1}_{p,q}~(J')$, $e_{J'} \in V^{0}_{n-1}$. (Note that $q \notin J$). Then to complete the proof it remains to show that the claim is valid for any such $J$ that does not belong to $\mathcal{I}_{n}$. 
\\

To implement the outlined strategy, we will consider several cases. We may assume that $q\neq n+1$. Indeed, the only index set $I$ such that $ [1...n+1] \setminus I = \{n+1\}$ is $I = [1...n]$. Then the bounded ratio $\frac{[1...n][2...n-1,n+1,2n]}{[1...n-1,n+1][2...n, 2n]}$ allows one to express $w_I$ via $w_{I_{1}}, w_{I_{2}},w_{I_{3}}$, where for each of the  three index sets $I_{1,2,3}$ present in the ratio 
%we observe that $w_{I}=w_{I_{1}}+w_{I_{2}}-w_{I_{3}}$ (mod $V_{n}$), 
%such that for $I_{1}$, $I_{2}$, $I_{3}$ 
there is a choice of $q\neq n+1$.

We start with the case when $n+1 \in I$. Then we can chose $p=n+1$, and some $q \in [1...n]$.
First, let $q\neq n$. Then the first $n-2$ elements $J$ of the set $\{{p}~\cup~\mathbf{S}^{-1}_{p,q}~(J')$, $J' \in \mathcal{I}_{n-1}\}$ are $J_i=\{i\} \cup [n+1,n+3...2n] (i=1,...,\hat{q},...,n-1)\}$. %For all such $J$ there is $w_{J} = e_{J}~(\text{mod} ~V^{0}_{n})$. Indeed, $e_{J} = e_{\tilde{J}}~(\text{mod} ~\overline{V}_{n}/V^{0}_{n})$ where $\tilde{J}=(J\setminus\{i\})\cup\{i+1\}$, which comes from 
Consider the bounded ratio $$\frac{[i+1,n+1,n+3...2n][i,n+2,n+3...2n]}{[i,n+1,n+3...2n][i+1,n+2,n+3...2n]}$$ and the corresponding vector $w_i\in \overline{V}_n$.
Then $e_{J_i}=e_{J_{i+1}} + w_i \ (\text{mod}\ V^{0}_{n})$. Continuing in the same way and observing that $J_n\in \mathcal{I}_{n}$, we conclude that there is $w_{J} = e_{J}~(\text{mod} ~V^{0}_{n})$.

%Thus, $e_{J} = %e_{\tilde{J}}=...=e_{[n,n+1,n+3...2n]}=0~%%(\text{mod} ~\overline{V}_{n}/V^{0}_{n})$. 
The remaining $n$ elements $J$ in $\{{p}~\cup~\mathbf{S}^{-1}_{p,q}~(J')$, $J' \in \mathcal{I}_{n-1}\}$ are 
$$\{\{[n+1,n+2,n+3...2n], [n,n+1,n+3...2n]\} \cup \{[n,n+1,n+2,...2n]\setminus{\{j\}} ,(j=n+3,...,2n)\} \}$$ 
and all of them are in $\mathcal{I}_{n}$.

Next, let $q = n$. Then the first $n-1$ elements $J$ in the set $\{{p}~\cup~\mathbf{S}^{-1}_{p,q}~(J')$, $J' \in \mathcal{I}_{n-1}\}$ -  $\{i\} \cup [n+1,n+3...2n]\ (i=1,...,n-1)$ - have been dealt with above, and 
%For all such $J$ there is $w_{J} = e_{J}~(\text{mod} ~V^{0}_{n})$, similar to the case $(a.1)$. 
the $n$-th element is $J=[n+1,n+2,...,n]\in \mathcal{I}_{n}$.
%from the set 
%$\{{p}~\cup~\mathbf{S}^{-1}_{p,q}~(J')$, $J' %\in \mathcal{I}_{n-1}\}$ is $[n+1,n+2,...,n]$ is in $\mathcal{I}_{n}$. 
The remaining $n-2$ elements in
$\{{p}~\cup~\mathbf{S}^{-1}_{p,q}~(J')$, $J' \in \mathcal{I}_{n-1}\}$ are $J_j=[n-1,n+1,n+2,...,2n] \setminus \{j\}\ (j=n+3,...,2n) \}$.
%there is such $w_{J} = e_{J}~(\text{mod} ~V^{0}_{n})$. Indeed, $e_{J} = e_{J'}~(\text{mod} \overline{V}_{n}/V^{0}_{n})$ where $J'=(J\setminus\{j-1\})\cup\{j\}$, which comes 
Consider the bounded ratio $$\frac{[n,n+1,n+2...\hat{j}...2n][n-1,n+1,n+2...\hat{(j-1)}...2n]}{[n-1,n+1,n+2...\hat{j}...2n][n,n+1,n+2...\hat{(j-1)}...2n]}$$
and the corresponding vector $w_i\in \overline{V}_n$.
%Thus, $e_{J} = e_{J'}=...=e_{[n-%1,n+2,n+3...2n]}=0~(\text{mod} %~\overline{V}_{n}/V^{0}_{n})$  
Then $e_{J_j}=e_{J_{j-1}} + w_j \ (\text{mod}\ V^{0}_{n})$. 
Continuing in the same way and observing that $J_{n+1}=[n-1,n+2,n+3...2n]$ is in 
$\mathcal{I}_{n}$, we conclude that there is $w_{J} = e_{J}~(\text{mod} ~V^{0}_{n})$.

Now let $n+1 \notin I$ and chose some $p \in [n+2...2n]$ and $q \in [1...n]$. 
If $q \neq n$, then $q \in [1,...,n-1]$ and $p \in [n+2,...,2n]$ and all elements in the set $\{{p}~\cup~\mathbf{S}^{-1}_{p,q}~(J')$, $J' \in \mathcal{I}_{n-1}\}$ are in $\mathcal{I}_{n}$.

If $q=n$, then the  first $n$ elements $J$ in the set $\{{p}~\cup~\mathbf{S}^{-1}_{p,q}~(J')$, $J' \in \mathcal{I}_{n-1}\}$ are $\{i\} \cup [n+2,...,2n], (i=1,...,n-1,n+1)$ which are all in $\mathcal{I}_{n}$. The remaining $(n-2)$ elements in $\{{p}~\cup~\mathbf{S}^{-1}_{p,q}~(J')$, $J' \in \mathcal{I}_{n-1}\}$ are $[n-1,n+1,n+2,...,2n] \setminus\{j\}\  (j=n+2,...,2n)$ which were  already treated above. 

%We obtained that $\mbox{dim}(\overline{\rm %\mathcal{C}}_{n}) \geq \binom{2n}{n}-2n$. 
This completes the proof by induction. 
\end{proof}

\subsection{Maximal independent subsets of primitive vectors 
%$\mathcal{P}_{n}$
} 
First, we consider linear dependencies between primitive vectors $\Vec{v}_{i,j,\Delta}$. 
%and show that dimensions of $\overline{\rm \mathcal{C}}_{n}$ and $\mathcal{C}_{n}$ are equal to $\binom{2n}{n}-2n$ for any $n$.

\begin{lem}\label{lm:ind}
(1) If $j-i=n$ and elements of $\Delta$ form an interval $[i+2,j-1]$, then $\Vec{v}_{i,j,\Delta}$ is linearly independent with the remaining  primitive  vectors.\\

(2) %If there is a dot without a cross on the circle to the left or to the right of a cross, then this primitive  vector satisfies the following linear equations with other primitive  vectors as it is described on the Fig. 5: 
%Take difference of the initial primitive  vector and a primitive  vector obtained by shifting a cross to it's empty adjacent dot on the circle. This difference is equal to the difference of other two vectors, there the the shifted pair from the first two vectors became a pair of circles, and a pair of circles from the first two vectors became a new shifted crosses in the same orientation as it was in the first pair, and the rest dots and circles are the same (see Fig.5). Note, that there are two choices of the pair of circles, so there are three equal differences. 
If $\Vec{v}_{i,j,\Delta}$ does not satisfy conditions in part $(1)$, then for any element $x\in\Delta$ such that $(x+1$ mod $2n)\notin\Delta$, $(x+1$ mod $2n) \neq i$, $(x+1$ mod $2n) \neq j$ the following relations hold%(see an example in Fig. \ref{fig:general}):
$$\Vec{v}_{i,j,\Delta}-\Vec{v}_{i,j,(\Delta \setminus x)\cup(x+1)} = \Vec{v}_{\min\{x,j\},\max\{x,j\},(\Delta \setminus x)\cup(i)}-\Vec{v}_{\min\{x,j\},\max\{x,j\},(\Delta \setminus x)\cup(i+1)}=$$
\begin{equation}\label{eq:2-2}
=\Vec{v}_{\min\{x,i\},\max\{x,i\},(\Delta \setminus x)\cup(j)} - \Vec{v}_{\min\{x,i\},\max\{x,i\},(\Delta \setminus x)\cup(j+1)}
\end{equation}
\end{lem}

% \begin{figure}[htb]
% \centering
% \includegraphics[width=5.1in]{f3.jpg}
% \caption{\\}
% \label{fig:general}
% \end{figure}
% (Above, blue arrows indicate how the second diagram in each difference is obtained from the first.)
\begin{proof} 
(1) $\Vec{v}_{i,j,\Delta}$ has a nonzero coordinate corresponding 
%to an axis associated with 
Pl\"ucker coordinate $[i+1,j]$. Clearly, for all other primitive  vectors the same coordinate is zero. Thus, $\Vec{v}_{i,j,\Delta}$ is linearly independent with the remaining primitive  vectors.
\\
(2) Relations \eqref{eq:2-2} follow in a straightforward way from Definition 2.3 of the primitive  vectors. %Each vector has four nonzero coordinates, and each relation is just a different grouping of eight  nonzero coordinates. 
\end{proof}
Note that relations \eqref{eq:2-2} remain valid if $x+1$ is replaced by $x-1$. 
%The proof follows from multiplication of each term in the relations (\ref{eq:2-2}) by $-1$.
Note also that sometimes
%primitive vectors have several different ways to shift one cross to the adjacent unmarked dot. Such 
the same primitive vector can be involved in several different linear relations of type (\ref{eq:2-2}). For $n=4$, a complete list of all relations consists of two relations of type (\ref{eq:2-2}) and all their cyclic shifts and reflections (Theorem 2.6).
$$\Vec{v}_{1, 3,\{5,6\}}-\Vec{v}_{1,3,\{5,7\}} = \Vec{v}_{1, 6,\{3,5\}}-\Vec{v}_{1,6,\{4,5\}} = \Vec{v}_{3, 6,\{1,5\}}-\Vec{v}_{3,6,\{2,5\}}$$
$$\Vec{v}_{2, 4,\{1,6\}}-\Vec{v}_{2,4,\{6,8\}} = \Vec{v}_{4, 8,\{3,6\}}-\Vec{v}_{4,8,\{2,6\}} = \Vec{v}_{2, 8,\{5,6\}}-\Vec{v}_{2,8,\{4,6\}}$$

The proposition below provides an explicit  description for a basis made of primitive vectors for the smallest linear subspace containing the cone $\mathcal{C}_{n}$.

%For each Pl\"ucker coordinate given by an index subset $I$ we associate a diagram consisting of a circle $c_{I}$ with $2n$ dots enumerated clockwise, where the dots labeled by elements from $I$ are marked by crosses.

To visualize the set of primitive vectors that form a basis, we associate to each  $\Vec{v}_{i,j,\Delta}$  a diagram $C_{i,j,\Delta}$ consisting of a circle  with $2n$ dots enumerated clockwise, where $(n-2)$ dots labeled by elements from $\Delta$ are marked by crosses, and two pairs of dots corresponding to pairs of $i,i+1$ and $j,j+1$ are circled (see an example below of the diagram associated to a primitive ratio
$\frac{[1238][1358]}{[1258][1368]}$). \begin{figure}[htb]
\centering
\includegraphics[width=2.in]{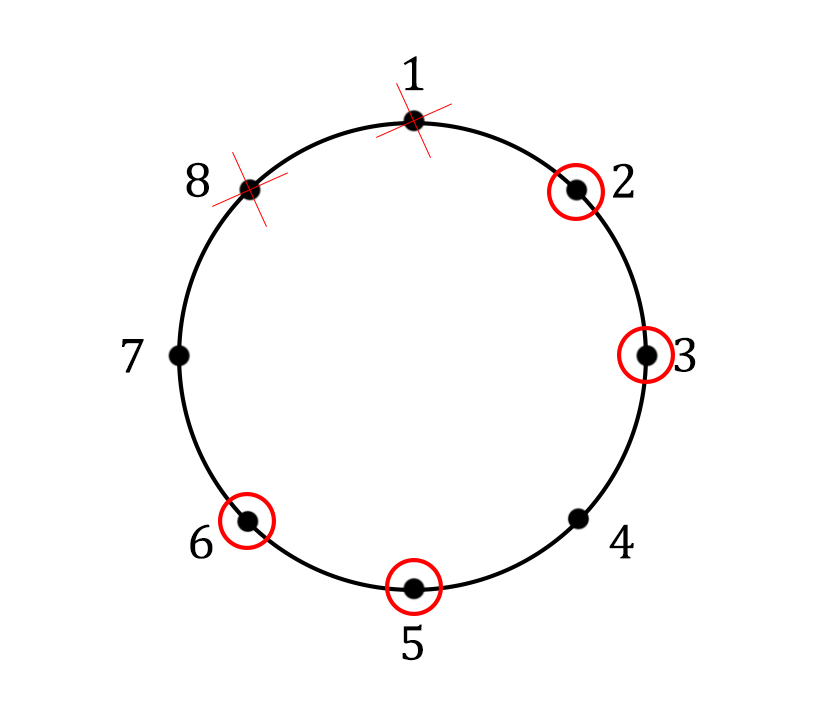}
%\caption{\\}
\label{fig:general}
\end{figure}

Given a diagram $C_{i,j,\Delta}$ (equivalently, $\Vec{v}_{i,j,\Delta}$), we define three arcs (equivalently, intervals)  $Arc_{1}=[1,i-1], Arc_{2}=[i+2,j-1], Arc_{3}=[j+2,2n]$. If $j=2n$ we define $Arc_{1}$ to be $[2...i-1]$ and $Arc_{3}$ to be empty. If $i=1$, $Arc_1$ is empty.
%Values of $i$ could be in the range $[1...2n-2]$. 
%In what follows we will sometimes abuse notations by using $c_I$ to mean both the diagram and the corresponding Pl\"ucker coordinate.

\begin{figure}[ht]
\begin{center}
\includegraphics[width=2.0in]{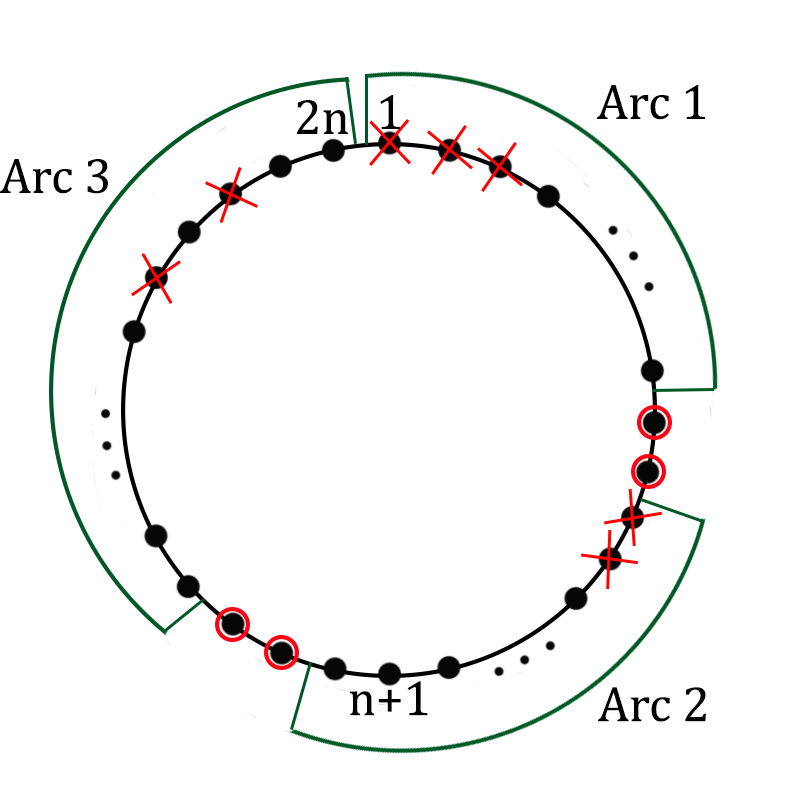}
\end{center}
\caption{}
\label{f18}
\end{figure}

Consider a subset $B_{n}$ of primitive vectors described by diagrams $C_{i,j,\Delta}$ that satisfy the following conditions: 
\begin{enumerate}
\item  If nonempty, $Arc_{1}\cup \Delta$ is a subinterval of $Arc_{1}$ that contains the initial point of $Arc_{1}$.
%If there are any crosses on the $Arc_{1}$, they all are consecutive and start at one.
\item If nonempty, $Arc_{1}\cup \Delta$ is a subinterval of $Arc_{2}$ that contains the initial point of $Arc_{2}$.
%If there are any crosses on the $Arc_{2}$, they all are consecutive and start at $i+2$.
%\item  The $Arc_{3}$, all crosses are in arbitrary positions.
\end{enumerate}
There are no restrictions on $Arc_{3}\cup \Delta$.
An example of an element from ${B}_{n}$ is in Fig.\ref{f18}.

\begin{prop}
 The set $B_{n}$ is a basis of the subspace of $V_n$ generated by vectors from $\overline{\mathcal{C}}_{n}$. 
\end{prop}
%\begin{remark} 

%\end{remark}   

\begin{pf}
By Proposition~\ref{th:dim}, it suffices to show that $\vert B_{n}\vert=\binom{2n}{n}-2n$ and that all primitive  vectors can be expressed as linear combinations of vectors from $B_{n}$.  

The number of elements in the set $B_{n}$ is given by (\ref{eq:Bn}) below\footnote{D. Galvin suggested an alternative counting method by establishing a bijection with the set of {\em crossing} $n$-element subsets of the $2n$-element set, whose  size is $\binom{2n}{n}-2n$.}. To see this, first note that $B_{n}$ contains $(n-1)^2$ elements such that $j=2n$. Indeed, it is not hard to see that if the set $Arc_{1}\cup \Delta$  is nonempty, its largest element of $Arc_{1}\cup \Delta$ is less than $n-1$ and, 
if  $Arc_{2}\setminus \Delta$  is nonempty, its smallest element is greater  than $n+1$ and less than $2n-1$. These gives us $(n-1)^2$ choices, each of which determines $\Vec{v}_{i,2n,\Delta}\in B_n$ uniquely.
 In the second term of \eqref{eq:Bn} , $k$ represents the size of $Arc_{3}\cup \Delta$ and determines the range of possible values for $j$ in the internal sum.  For a fixed $k$, the  binomial coefficient counts the number of possibilities for  $Arc_{3}\cup \Delta$. Furthermore, the largest element of $Arc_{1}\cup \Delta$ cannot exceed $n-k-2$, which gives $(n-k-1)$ choices and the smallest element of $Arc_{2}\setminus \Delta$ can be between $n-k+1$ and $j-1$, which gives $(j-n+k)$ choices. This justifies the formula below:
 
 %The number of possible positions of crosses on the $Arc_{1}$ is $n-k-1$, the number of possible positions of unmarked dots on the $Arc_{1}$ is $j-n-k$.
\begin{equation}\label{eq:Bn}|B_{n}| = (n-1)^2 + \sum_{k=0}^{n-2}  \sum_{j=n-k+1}^{2n-k-1} \binom{2n-j-1}{k}(n-k-1)(j-n+k)
\end{equation}
or, relabeling $j-(n-k)$ as $j$,
$$ |B_{n}| = (n-1)^2 + \sum_{k=0}^{n-2}  \sum_{j=1}^{n-1} \binom{n-j+k-1}{k}(n-k-1)j$$
%We assume it is equal to $\binom{2n}{n}-2n$. 
To prove that $\vert B_{n}\vert=\binom{2n}{n}-2n$, we argue by induction and compute
%,  to show that for $n+1$ the following holds:
%$$ |B_{n+1}| = n^2 + \sum_{k=0}^{n-1}  %\sum_{j=1}^{n} \binom{n-j+k}{k}(n-k)j = %\binom{2n+2}{n+1}-2(n+1)$$
%it is enough to show that $|B_{n+1}|-|B_{n}| %= \binom{2n+2}{n+1}-2(n+1) - \binom{2n}%{n}+2n$.
\\
$|B_{n+1}|-|B_{n}| =$
$$= n^2 + \sum_{k=0}^{n-1}  \sum_{j=1}^{n} \binom{n-j+k}{k}(n-k)j - (n-1)^2 - \sum_{k=0}^{n-2}  \sum_{j=1}^{n-1} \binom{n-j+k-1}{k}(n-k-1)j = $$
$$ = (2n-1) + \sum_{k=0}^{n-1}  \sum_{j=1}^{n} \binom{n-j+k}{k}(n-k)j - \sum_{k=0}^{n-2}  \sum_{j=2}^{n} \binom{n-j+k}{k}(n-k-1)(j-1) = $$
$$ = (2n-1)+ \sum_{k=0}^{n-2}  \sum_{j=2}^{n} \binom{n-j+k}{k}(n-k)j + \sum_{j=1}^{n} \binom{2n-j-1}{n-1}j + \sum_{k=0}^{n-2}\binom{n+k-1}{k}(n-k)-$$ 
$$-\sum_{k=0}^{n-2}  \sum_{j=2}^{n} \binom{n-j+k}{k}(n-k-1)j +\sum_{k=0}^{n-2}  \sum_{j=2}^{n} \binom{n-j+k}{k}(n-k-1) = $$
$$ = (2n-1) + \sum_{k=0}^{n-2}  \sum_{j=2}^{n} \binom{n-j+k}{k}j + \sum_{j=1}^{n} \binom{2n-j-1}{n-1}j + \sum_{k=0}^{n-2}\binom{n+k-1}{k}(n-k) + $$
$$ + \sum_{k=0}^{n-2}  \sum_{j=2}^{n} \binom{n-j+k}{k}(n-k+1) = (2n-1) + S_{1} + S_{2}+ S_{3}+ S_{4}.$$
Using well known identities for binomial coefficients, one can show that
\\
$S_{1} = \binom{2n-1}{n-1} + \binom{2n-2}{n-1} -n -1,
S_{2} = \binom{2n}{n+1},
S_{3} = \binom{2n-1}{n+1}+\binom{2n-2}{n},
S_{4} =  \binom{2n-1}{n}-n \ .$
Thus, $(2n-1) + S_{1} + S_{2}+ S_{3}+ S_{4} = -2 - \binom{2n}{n} + \binom{2n+2}{n+1}.$

Now, we express all primitive  vectors as linear combinations of vectors from the set $B_{n}$. Observe, that all primitive vectors containing a nonzero coordinate corresponding to a cyclically dense Pl\"ucker coordinate are already in the set $B_{n}$. All other primitive vectors are involved in at least one relation of type $v_{1}-v_{2}=v_{3}-v_{4}$ of the kind described in Lemma \ref{lm:ind}

Let us fix a primitive vector $v$. If it does not satisfy conditions (1) and (2), we will express express it as  a linear combination of vectors from $B_{n}$ following the steps described below.
\\
%If the initial vector satisfies (1), then nothing needs to be done with the $Arc_{1}$. Otherwise,
If $v$ does not satisfy the condition (1), we select the first cross that follows the first unmarked dot on $Arc_{1}$ in a clockwise order. We shift this cross counterclockwise to the unmarked dot, as illustrated in Fig.\ref{f19}. 
\begin{figure}[htb]
\centering
\includegraphics[width=5.1in]{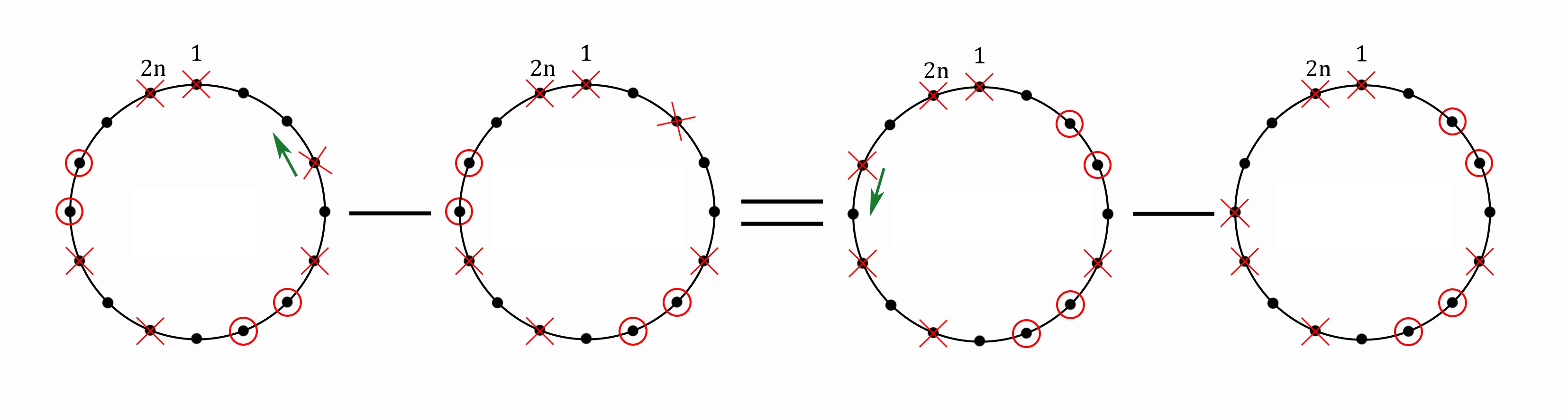}
\caption{ }
\label{f19}
\end{figure}
Here we use a relation of type $v_{1}-v_{2}=v_{3}-v_{4}$ that involves $v$, the shift of $v$ that is lexicographically closer to a vector from $B_{n}$
%, and the other two vectors which have a longer $Arc_{3}$. Observe that the shifted vector is lexicographically closer to a vector from $B_{n}$, 
and two more primitive vectors that satisfy (1). 
%and their 
($Arc_{3}$ for both of them is longer $Arc_{3}$ for $v$.) 
It follows that $v$  is expressed as a linear combination of primitive vectors with $Arc_{1}$ satisfying condition (1) for $B_n$. For any of these vectors, a similar procedure can be applied to $Arc_{2}$. This time, we use 
%
%Now, repeat the same procedure with the $Arc_{2}$. Starting from the first dot of the $Arc_{2}$ we move clockwise to the first cross going after unmarked dot on the $Arc_{2}$. Then we shift this cross counterclockwise to the unmarked dot. There is
a relation of type $v_{1}-v_{2}=v_{3}-v_{4}$ which contains a vector that satisfies (1) but not (2), its shift that is lexicographically closer to a vector from $B_{n}$ and two other coordinates 
%whose $Arc_{3}$ are longer. Observe that the last two coordinates are from the set
that are in $B_{n}$. We conclude that after finitely many of such  counterclockwise shifts any primitive vector can be expressed as a linear combination of elements from $B_{n}$. By Theorem \ref{th:dim}, the dimension  of the span of all primitive vectors is $\binom{2n}{n}-2n$, and therefore the same is true for the span of vectors in $B_{n}$, which means that $B_{n}$ is a basis. $\diamond$
\end{pf}

\subsection{System of inequalities for $\overline{\rm \mathcal{C}}_{n}$}

In this subsection we describe the method of finding the system of linear inequalities that defines $\overline{\rm \mathcal{C}}_{n}$. Then we specialize to $n=4$ to show that Conjecture~\ref{conj:mainconj} part $(1)$ does not hold. 

Recall that a finite system of linear inequalities is reduced if there is none of the  inequalities in the system which can be represented as a non-negative linear combination of other inequalities in it.
Each inequality of a reduced system corresponds to a facet of the cone of solutions to this system. 
%Indeed, linear subspace of codimension one given by boundary of the half-space defined by the inequality intersects the cone of solutions to this system by facet of codimension one. 
There is a one-to-one correspondence between facets of a cone and subsets of its extreme rays with the following three properties: 
\begin{enumerate}
\item  Linear subspace generated by vectors in the subset has a co-dimension one in the cone.
\item All extreme rays not in the subset belong to a half-space bounded by the hyperplane that contains the given facet. 
\item Such subset is maximal by inclusion.
\end{enumerate}
%Note, that the second property assures that the generated hyperplane is facet of the cone, not its diagonal section. 
We call subsets of extreme rays with these three properties  {\em facet sets}. Subsets of extreme rays which are complements to facet sets will be  called {\em outer sets}. The linear inequality defined by an outer set $U$ (equivalently, by the corresponding facet set) is 
$(\Vec{n},\Vec{x})<0$, where $\Vec{n}$ is a vector normal to the subspace spanned by the facet set that satisfies $(\Vec{n},\Vec{f})<0$ for any $\Vec{f}\in U$.

%For each outer set $U$ we describe the corresponding inequality in the system which defines the cone spanned by all primitive vectors. Let $\Vec{n} \in V$ be such that $(\Vec{n},\Vec{u})=0$ for any $\Vec{u} \in U$ and $(\Vec{n},\Vec{f})<0$ for any $\Vec{f}$ in the corresponding facet set. Then, the inequality from the system defining $ \mathcal{C}_{3}$ is $(\Vec{n},\Vec{x})<0$, $\Vec{x} \in V$

Using the three properties listed above, one can find all facets of $\overline{\rm \mathcal{C}}_{n}$ by exhaustive search.  
%All primitive vectors can not belong to a single facet. 
For each primitive vector we are searching for all outer sets containing this primitive vector. There are two possible cases: 

$1.$ Initial primitive vector is one of the $2n$ vectors containing a coordinate corresponding to a cyclically dense Pl\"ucker variable. Then, the rest of the vectors yield a facet set, since they span a subspace of dimension one less then the dimension of the cone, and the properties (2) and (3) hold because the outer set consists of a single vector. 

$2.$ Initial primitive vector is not as in case $1$. Then, it is involved in at least one relation $v_{1}-v_{2}=v_{3}-v_{4}=v_{5}-v_{6}$ of type \eqref{eq:2-2} 
%(for the case $n=3$ it is exactly one relation of the given type). 
Without loss of generality, we can assume that the initial vector is $v_{1}$. An outer set containing vector $v_{1}$ must contain either vector $v_{2}$, or $v_{3}$ and $v_{5}$. Otherwise, $v_2$ is in a facet set and then either vector $v_1$ is a linear combination of three vectors from a facet set, or one the following two equations holds $(v_1,\Vec{n})=-(v_4,\Vec{n})$,
$(v_1,\Vec{n})=-(v_6,\Vec{n})$, where $\Vec{n}$ is the normal to the hyperplane containing the facet. The first case contradicts the assumption that $v_{1}$ is in outer set, the second case violates property $(2)$. Thus we find another primitive vector in the outer set containing $v_1$. We then apply the same process to this new vector using other linear dependencies \eqref{eq:2-2} until the whole outer set is found. Note that, since $v_1$ can be involved in several relations \eqref{eq:2-2}, the procedure above may have to be repeated.
%In other words, given that some vector is in an outer set, in all the relations $v_{1}-v_{2}=v_{3}-v_{4}=v_{5}-v_{6}$ involving this vector either another vector in the same part of the double equation must belong the outer set as well, or other two vectors in two other parts of the double equation with the same sign as the initial deleted vector belong to the outer set.\\

\begin{figure}[htb]
\centering
\includegraphics[width=5.5in]{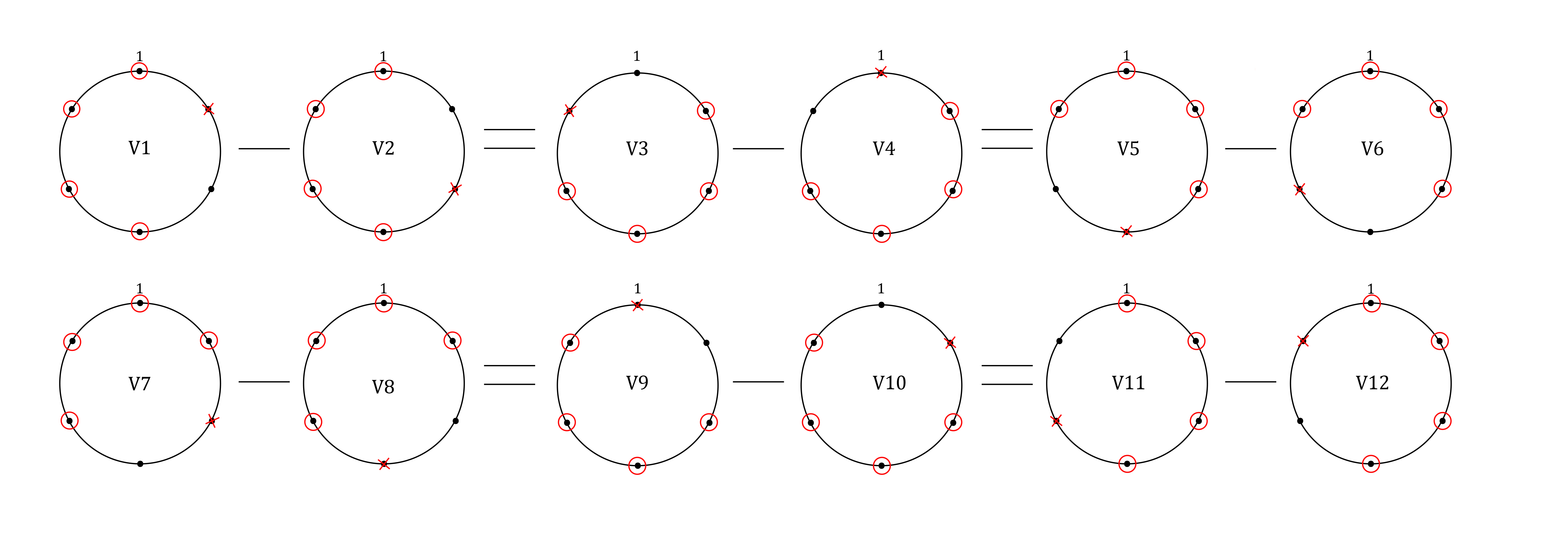}
\caption{}
\label{f20}
\end{figure}

Let us  illustrate this method in the smallest case $n=3$.
There are $18$ primitive vectors. We label them in such a way that 
%{\color{blue}
%\[
%\{v_1, \ldots, v_{18}\}
%= list\ explicitly\ using \ notations\  v_{i,j,\Delta} !!!
%\]
%}
the last six involve cyclically dense Plücker coordinates and thus are independent with the rest of primitive vectors by Lemma~\ref{lm:ind} and so, each one of defines an outer set. Remaining 12 vectors are subject to relations depicted in Fig.\ref{f20}.

The cone $\overline{\rm \mathcal{C}}_{3}$ the case $n=3$ has $16$ facets. The corresponding $16$ outer sets can be easily found using 
Fig.\ref{f20} : 
$$
\{v_{13}\}, \{v_{14}\}, \{v_{15}\}, \{v_{16}\}, \{v_{17}\}, \{v_{18}\}, \{v_{1}, v_{2}\}, \{v_{3}, v_{4}\}, \{v_{5}, v_{6}\}, \{v_{7}, v_{8}\}, $$
$$\{v_{9}, v_{10}\}, \{v_{11}, v_{12}\}, \{v_{1}, v_{3}, v_{5}\}, \{v_{2}, v_{4}, v_{6}\}, \{v_{7}, v_{9}, v_{11}\},\{v_{8}, v_{10}, v_{12}\}\ .
$$

Recall, that  Proposition 2.8 implies that the cones $\overline{\rm \mathcal{C}}_{3}$ and $\mathcal{C}_{3}$ coincide. We reaffirm this fact by verifying that the system of linear inequalities for $\overline{\rm \mathcal{C}}_{3}$ defined by $16$ outer sets above is identical to
the system of linear inequalities for $\
{\rm \mathcal{C}}_{3}$ obtained
following the method outlined at the end of Section 3. 

In contrast, when we applied the two methods above to compute reduced systems of linear inequalities for $\overline{\rm \mathcal{C}}_{4}$ and $\mathcal{C}_{4}$, it turned out that the latter is a strict subsystems of the former. Therefore, 
$\mathcal{C}_{4}   \supsetneq \overline{\rm \mathcal{C}}_{4}$ which contradicts 
part (1) of Conjecture~\ref{conj:mainconj}.

%Similarly, one can find outer subsets corresponding to $\overline{\rm \mathcal{C}}_{4}$. The only difference is that some primitive vectors belong to several relations. In order to find an outer sets containing such vectors we have to repeat the procedure described above with all relations involving the initial vector. We verified that in each case, the subset obtained by the above process which must belong to an outer set was an outer set itself. Thus, we found the system of inequalities of $\overline{\rm \mathcal{C}}_{4}$, which is a reduced system.

%Now we make a key observation, that the system of linear inequalities which corresponds to the obtained list of facets of $\overline{\rm \mathcal{C}}_{4}$ is not the same as the system of inequalities of $\mathcal{C}_{4}$. In fact, the system of inequalities which corresponds to facets of $\mathcal{C}_{4}$ is a strict subsystem of the reduced system of $\overline{\rm \mathcal{C}}_{4}$. Thus, we conclude that part $(1)$ of Conjecture~\ref{conj:mainconj} does not hold. Instead, $\mathcal{C}_{4}   \supsetneq \overline{\rm \mathcal{C}}_{4}$. 

Furthermore, we found several bounded ratios which are extreme rays of $ \mathcal{C}_{4}$ and do not belong to $\overline{\rm \mathcal{C}}_{4}$. They are described in the next section.

\section{New inequalities}
%In this section we describe two numerical methods which provide examples of extreme rays in $\mathcal{C}_{4}$, which are not primitive vectors.
%
Since finding the full list of extreme rays for $\mathcal{C}_{4}$ proved to be not feasible computationally, we needed to develop a method that allows you to find at least some of the extreme rays that are not primitive vectors and five rays to genuinely new bounded ratios.

We start with the following lemma that provides as useful tool to check if a vector is an extreme vector of a cone.
\begin{lem}
Consider a vector $v$ that satisfies a reduced system of homogeneous linear inequalities 
$\{Ax \leq 0\}$. Assume further that the corank of the submatrix $A'$ of $A$ such that $A'v=0$ is equal to $1$.  Then $v$ is an extreme vector of the cone defined by $A$.
\end{lem}

\begin{pf}
Assume that $v$ is not an extreme vector of the cone $\{Ax \leq 0\}$. Then $v$ has to be a linear combination with positive coefficients of at least two  vectors in the cone non-collinear to $v$.  Then there exists a 2-dimensional disk in the linear subspace spanned by two vectors in this linear combination, such that it is contained in the cone and is centered at $\alpha v$. For each hyperplane which is annihilated by $v$, there are two possibilities: either the whole disc is contained in the hyperplane or there exists a point $\alpha v + \epsilon$ in the disc which does not belong to the hyperplane. In the latter case, $\alpha v - \epsilon$ does not belong to the hyperplane as well. Then, two points $\alpha v + \epsilon$ and $\alpha v - \epsilon$ are in two different half spaces with respect to the hyperplane. This contradicts the fact that the whole disc belongs to the cone. Thus, the former possibility is the case for all hyperplanes  annihilated by  $v$ and therefore the whole  2-dimensional disc belongs to the kernel of $A'$ of corank 1, contradiction. So, $v$ is an extreme vector.$\qed$
\end{pf}
%The lemma above is a practical tool to verify if a vector is an extreme vector. 

To discover new extreme vectors of $\mathcal{C}_{4}$, 
we first pick an inequality $l x\leq0$ from the system $\{Ax\leq0\}$ defining $\overline{\rm \mathcal{C}}_{4}$ which is missing in the system defining $\mathcal{C}_{4}$. Then we consider a new system $\{A'x\leq0\}$ obtained from $\{Ax\leq0\}$ by replacing $l x\leq0$ with $-l x\leq0$. Any vector satisfying this system is inside $\mathcal{C}_{4}$ but outside $\overline{\rm \mathcal{C}}_{4}$. We call the cone of solutions to this system $C'$. To look for extreme rays of $C'$, we first pick an affine hyperplane that separates a bounded part of $C'$ from infinity. Then we minimize a random linear functional on that bounded part to get new extreme rays. We then verify if they are extreme rays for $\mathcal{C}_{4}$ using Lemma 5.1. 

Note that if a new extreme vector is found then all its cyclic shifts and reflections are extreme rays as well. It follows from Theorem 2.6.

Next, we recursively apply the approach above to a sequence of cones $\{K_{i}\}$, where $K_{1}=\overline{\rm \mathcal{C}}_{4}$ and  $K_{i+1}=K_{i}\cup E_{i}$, and $E_{i}$ is a set of new extreme rays obtained by applying our procedure  to $K_{i}  \subseteq \mathcal{C}_{4}$.  This method depends on whether we can find a facet of $K_{i+1}$ which corresponds to an inequality missing in the systems of defining inequalities for $\mathcal{C}_{4}$ and $K_{i}$. Unfortunately, after several steps this process becomes computationally unfeasible. This is the reason why the theorem below contains only a partial list of extreme vectors for for $\mathcal{C}_{4}$
%Such facets were possible to find while the new extreme rays for $\mathcal{C}_{4}$ were of particular type. The first eight new extreme rays can be expressed as a linear combination of primitive vectors with only one minus sign. The ninth new extreme ray has at least two minus signs, and computations become not feasible.}

A list of several new extreme rays none of which are related by the symmetries of Theorem \ref{thm:symm} is given in
\begin{thm}\label{t:main}
All ratios below are subtraction-free in terms of positive face weights, and thus are bounded by 1. Moreover, all of them can not be factored into a product of primitive ratios and thus provide explicit counterexamples to Conjecture~\ref{conj:mainconj}~part~(1).
\begin{enumerate}
%\Large 
\item $\frac{[1 3 6 8][1 4 5 8][1 4 6 7][2 3 5 8][2 3 6 7]}{[1 3 5 8][1 3 6 7][1 4 6 8][2 3 6 8][2 4 5 7]}\leq 1$
\vspace{0.3cm} 
\item $ \frac{[1 2 4 6][1 2 5 7][1 3 5 6][2 5 6 8][3 4 7 8]}{[1 2 5 6][1 3 4 6][1 3 5 7][2 4 6 8][2 5 7 8]}\leq 1$
\vspace{0.3cm}
\item $ \frac{[1 2 4 7][1 2 5 6][1 3 4 6][2 5 7 8][3 5 6 8]}{[1 2 4 6][1 2 5 7][1 3 5 6][2 5 6 8][3 4 7 8]}\leq 1$
\vspace{0.3cm}
\item $ \frac{[1 3 5 8][1 4 5 7][1 4 6 8][2 3 6 7][2 4 5 8]}{[1 3 5 7][1 3 6 8][1 4 5 8][2 4 5 7][2 4 6 8]}\leq 1$
\vspace{0.3cm}
\item $ \frac{[1 3 4 8][1 3 6 7][1 4 5 7][1 4 6 8][2 3 4 7][2 3 5 8][2 4 6 7]}{[1 3 4 7][1 3 5 7][1 4 5 8][1 4 6 7][2 3 4 8][2 3 6 7][2 4 6 8]}\leq 1$
\vspace{0.3cm}
\item $ \frac{[1 2 4 8][1 2 5 7][1 2 6 8][1 3 4 7][1 3 5 8][1 3 6 7][1 4 5 6][2 3 4 6][2 3 4 8][2 3 5 7][2 6 7 8][3 5 7 8][4 6 7 8]}{[1 2 4 7][1 2 5 8][1 2 6 7][1 3 4 6][1 3 4 8][1 3 6 8][1 4 6 7][2 3 4 7][2 3 5 6][2 3 5 8][2 5 7 8][3 6 7 8][4 5 7 8]}\leq 1$
\vspace{0.3cm}
\item $ \frac{[1 2 5 7][1 2 6 8][1 3 4 8][1 3 5 6][1 3 6 7][1 4 5 7][2 3 4 6][2 3 5 7][2 4 5 6][2 6 7 8][3 5 7 8][4 5 6 8][4 6 7 8]}{[1 2 5 8][1 2 6 7][1 3 4 6][1 3 6 8][1 4 5 6][1 4 6 7][2 3 4 7][2 3 5 6][2 4 5 7][2 5 7 8][3 5 6 8][3 6 7 8][4 5 7 8]}\leq 1$
\vspace{0.3cm}
\item $ \frac{[1 2 5 7][1 2 6 8][1 3 4 8][1 3 4 8][1 3 5 6][1 3 6 7][1 4 5 7][2 3 4 6][2 3 5 7][2 4 5 6][2 6 7 8][3 5 7 8][4 5 6 8]}{[1 2 4 8][1 2 6 7][1 3 4 6][1 3 5 7][1 3 6 8][1 4 5 6][1 4 6 7][2 3 4 8][2 3 5 6][2 4 5 7][2 5 7 8][3 5 6 8][3 6 7 8]}\cdot$
\\
$\frac{[4 6 7 8]}{[4 5 7 8]}\leq 1$
%\cdot$
%\\
%$\cdot \frac{[4 5 6 8][4 6 7 8]}{[3 6 7 8][4 5 7 8]}$
\vspace{0.3cm}
\item $ \frac{[1 2 4 7][1 2 5 7][1 2 6 8][1 3 4 8][1 3 4 8][1 3 5 6][1 3 5 6][1 3 6 7][1 4 5 7][2 3 4 6][2 3 5 7][2 4 5 6][2 5 6 8][2 6 7 8]}{[1 2 4 8][1 2 5 6][1 2 6 7][1 3 4 6][1 3 5 7][1 3 5 7][1 3 6 8][1 4 5 6][1 4 6 7][2 3 4 8][2 3 5 6][2 4 5 7][2 4 6 8][2 5 7 8]}\cdot$
\\
$\cdot\frac{[3 4 6 8][3 5 7 8][3 5 7 8][4 5 6 8][4 6 7 8]}{[3 4 7 8][3 5 6 8][3 5 6 8][3 6 7 8][4 5 7 8]}\leq 1$
\end{enumerate}
\end{thm}
\begin{proof}
The ratios listed above correspond to new extreme rays of the cone of bounded ratios $\mathcal{C}_{4}$, which do not belong to the cone of primitive ratios $\overline{\mathcal{C}}_{4}$. Moreover, a direct verification show that  all of them are subtraction-free in terms of face weights, thus they are bounded by 1. 
%These ratios disprove Conjecture~\ref{conj:mainconj}~part~(1). 
\end{proof}
Note, that ratios listed above give new determinantal inequalities for $M_{4}^{>0}$. For example, ratio $(3)$ is equivalent to 
$$\text{det}(A_{124|134})\text{det}(A_{12|12})\text{det}(A_{134|124})\text{det}(A_{2|3})\text{det}(A_{3|2}) \leq $$ 
$$\leq \text{det}(A_{124|124})\text{det}(A_{12|13})\text{det}(A_{13|12})\text{det}(A_{2|2})\text{det}(A_{34|34}) $$
\section{Concluding Remarks}
In order to explicitly describe $\mathcal{C}_{n}$ we need to find an intrinsic description of its extreme rays, that is to understand a `rule` which tells which ratios of minors correspond to extreme rays. The set of primitive ratios comes from short Pl\"ucker relations, which are closely related to mutations in a cluster algebra of $Gr(n,2n)$ \cite{scott2003grassmannians}. In this cluster algebra extended clusters that consist of Pl\"ucker coordinates are given by so called weekly-separated sets \cite{fomin2021introduction}. Thus, it is natural to associate a graph to each extreme ray, in which Pl\"ucker coordinates label the vertices and two vertices are connected by an edge iff Pl\"ucker coordinates associated to them are weekly-separated.We plan to investigate the structure of such graphs. An interesting initial observation is that for some extreme rays which are {\em not} related by symmetries of Theorem \ref{thm:symm}, corresponding associated graphs are nevertheless isomorphic (clearly, graphs associated to the extreme rays related by symmetries are isomorphic, since the property of being weekly-separated is preserved by symmetries).  \begin{figure}[htb]
\centering
\includegraphics[width=4.2in]{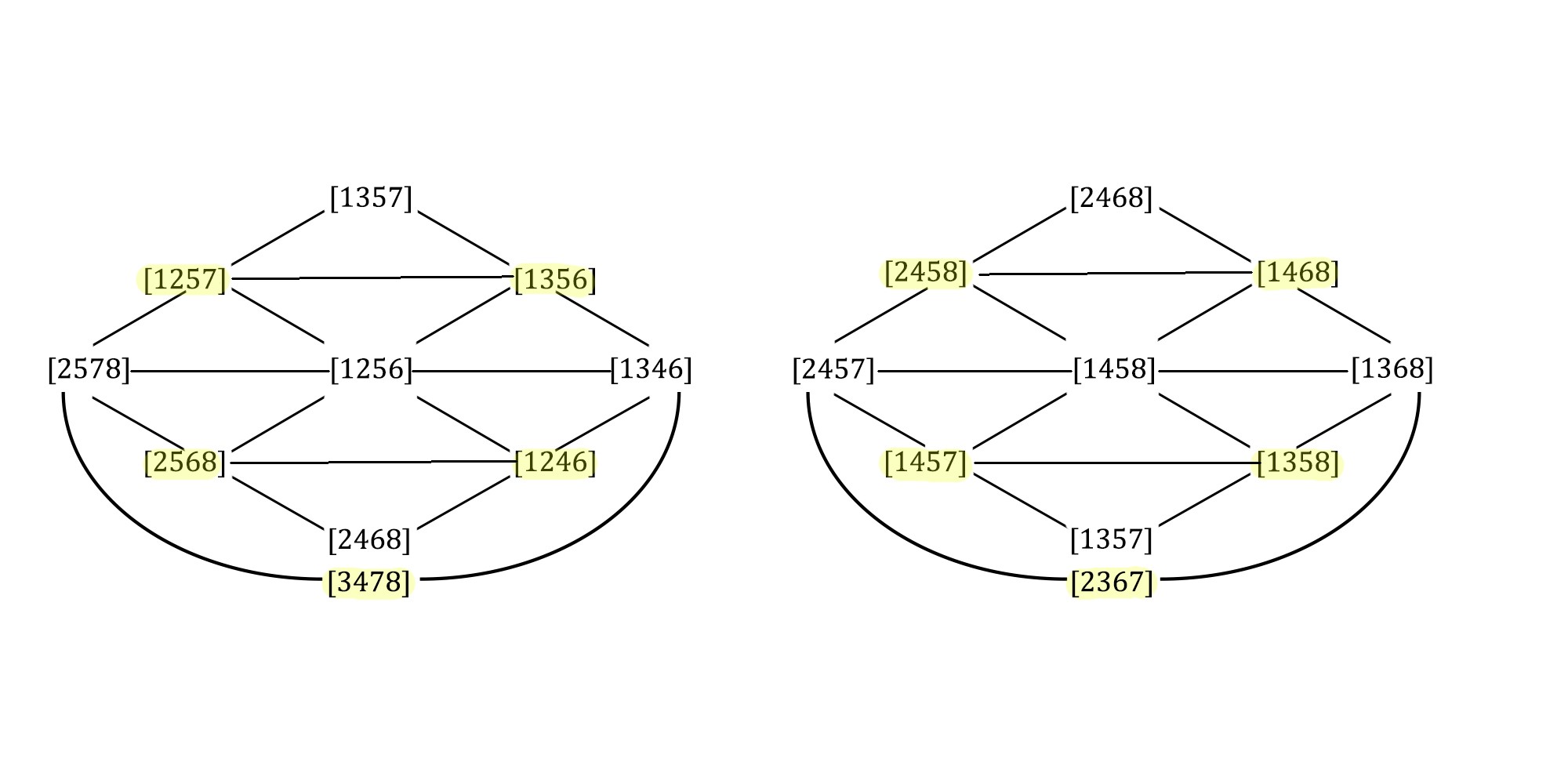}
\caption{}
\label{f21}
\end{figure} 
In Fig.\ref{f21}, we present two isomorphic graphs associated to the second and fourth extreme rays from the list of Theorem \ref{t:main}, where highlighted Pl\"ucker coordinates correspond to the numerator and white to the denominator of the corresponding ratio.
\section{Acknowledgements}
The authors would like to thank Komei Fukuda, David Galvin, Eric Jovinelly, Bohdan Kivva for inspiring and helpful conversations. This work was supported in part by the NSF grant DMS-2100785.

\bibliographystyle{plainnat}
\bibliography{m}
\end{document}